

\input amssym  

\newdimen\normalparindent

\iffalse  

\hsize=15truecm
\hoffset=.46truecm
\vsize=23.7truecm
\voffset=.46truecm

\normalparindent=24pt

\font\elevenrm=cmr10 at 11pt 
\font\eightrm=cmr8 
\font\sixrm=cmr6 

\font\eleveni=cmmi10 at 11pt 
\font\eighti=cmmi8
\font\sixi=cmmi6

\font\elevensy=cmsy10 at 11pt 
\font\eightsy=cmsy8
\font\sixsy=cmsy6

\font\elevenex=cmex10 at 11pt 

\font\elevenbf=cmbx10 at 11pt 
\font\eightbf=cmbx8
\font\sixbf=cmbx6

\font\eleventt=cmtt10 at 11pt 
\font\elevensl=cmsl10 at 11pt 
\font\elevenit=cmti10 at 11pt 

\textfont0=\elevenrm \scriptfont0=\eightrm \scriptscriptfont0=\sixrm
\def\rm{\fam0\elevenrm}
\textfont1=\eleveni \scriptfont1=\eighti \scriptscriptfont1=\sixi
 
\textfont2=\elevensy \scriptfont2=\eightsy \scriptscriptfont2=\sixsy

\textfont3=\elevenex \scriptfont3=\elevenex \scriptscriptfont3=\elevenex
\textfont\itfam=\elevenit
\def\it{\fam\itfam\elevenit}
\textfont\slfam=\elevensl
\def\sl{\fam\slfam\elevensl}
\textfont\bffam=\elevenbf \scriptfont\bffam=\eightbf
\scriptscriptfont\bffam=\sixbf
\def\bf{\fam\bffam\elevenbf}
\textfont\ttfam=\eleventt
\def\tt{\fam\ttfam\eleventt}

\skewchar\eleveni='177 \skewchar\eighti='177 \skewchar\sixi='177
\skewchar\elevensy='60 \skewchar\eightsy='60 \skewchar\sixsy='60

\font\sc=cmcsc10 at 11pt 

\font\cmcyr=cmcyr10 at 11pt
\font\cmcti=cmcti10 at 11pt
\font\cmccsc=cmccsc10 at 11pt

\smallskipamount=3.5pt plus 1pt minus 1pt
\medskipamount=7pt plus 2pt minus 2pt
\bigskipamount=14pt plus 2pt minus 2pt
\normalbaselineskip=14pt
\normallineskip=1pt
\normallineskiplimit=0pt
\jot=3.5pt

\normalbaselines
\rm

\else  

\hsize=15truecm
\hoffset=.46truecm
\vsize=23.7truecm
\voffset=.46truecm

\normalparindent=20pt

\font\sc=cmcsc10

\font\cmcyr=cmcyr10
\font\cmcti=cmcti10
\font\cmccsc=cmccsc10

\fi

\def\cyrchardefs{%
\chardef\yo=60
\chardef\Yo=62
\chardef\yu=64
\chardef\a=65
\chardef\b=66
\chardef\ts=67
\chardef\d=68
\chardef\ye=69
\chardef\f=70
\chardef\g=71
\chardef\kh=72
\chardef\i=73
\chardef\j=74
\chardef\k=75
\chardef\l=76
\chardef\m=77
\chardef\n=78
\chardef\o=79
\chardef\p=80
\chardef\ya=81
\chardef\r=82
\chardef\s=83
\chardef\t=84
\chardef\u=85
\chardef\zh=86
\chardef\v=87
\chardef\soft=88
\chardef\y=89
\chardef\z=90
\chardef\sh=91
\chardef\e=92
\chardef\shch=93
\chardef\ch=94
\chardef\hard=95
\chardef\Yu=96
\chardef\A=97
\chardef\B=98
\chardef\Ts=99
\chardef\D=100
\chardef\Ye=101
\chardef\F=102
\chardef\G=103
\chardef\H=104
\chardef\I=105
\chardef\J=106
\chardef\K=107
\chardef\L=108
\chardef\M=109
\chardef\N=110
\chardef\O=111
\chardef\P=112
\chardef\Ya=113
\chardef\R=114
\chardef\S=115
\chardef\T=116
\chardef\V=119
\chardef\Soft=120
\chardef\Y=121
\chardef\Z=122
\chardef\Sh=123
\chardef\E=124
\chardef\Shch=125
\chardef\Ch=126
\chardef\Hard=127
}
\def\cyr{\cmcyr\cyrchardefs}
\def\cycsc{\cmccsc\cyrchardefs}
\def\cyti{\cmcti\cyrchardefs}

\def\S{\mathhexbox278\thinspace}
\def\SS{\mathhexbox278\mathhexbox278\thinspace}

\def\square{\hbox to.77778em{%
\hfil\vrule\vbox to.675em{\hrule width.6em\vfil\hrule}\vrule\hfil}}

\long\def\acknowledgements#1\par{\medbreak\noindent{\it
    Acknowledgements\/}.\enspace #1\par\medbreak}
\def\definition#1\par{\medbreak\noindent{\bf Definition.}\enspace
  #1\par\medbreak}
\def\example#1\par{\medbreak\noindent{\bf Example.}\enspace
  #1\par\medbreak}
\long\def\remark#1\par{\medbreak\noindent{\it Remark\/}.\enspace
#1\par\medbreak}
\long\def\remarks#1\par{\medbreak\noindent{\it Remarks\/}.\enspace
#1\par\medbreak}
\def\exercise#1\par{\medbreak\noindent{\bf Exercise.}\enspace
#1\par\medbreak}
\def\notation#1\par{\medbreak\noindent{\bf Notation.}\enspace
#1\par\medbreak}
\def\proof{\noindent{\it Proof\/}.\enspace}
\def\endproof{\nobreak\hfill\quad\square\par\medbreak}

\def\lineover#1{{\offinterlineskip\mathchoice
{\setbox0=\hbox{$\displaystyle#1$}%
\vbox{\kern .33pt\hbox to\wd0{\kern 1pt\leaders\hrule height .33pt%
\hfill\kern 1pt}\kern 1pt\box0}}
{\setbox0=\hbox{$\textstyle#1$}%
\vbox{\kern .33pt\hbox to\wd0{\kern 1pt\leaders\hrule height .33pt%
\hfill\kern 1pt}\kern 1pt\box0}}
{\setbox0=\hbox{$\scriptstyle#1$}%
\vbox{\kern .25pt\hbox to\wd0{\kern .8pt\leaders\hrule height .25pt%
\hfill\kern .8pt}\kern .8pt\box0}}
{\setbox0=\hbox{$\scriptscriptstyle#1$}%
\vbox{\kern .2pt\hbox to\wd0{\kern .6pt\leaders\hrule height .2pt%
\hfill\kern .6pt}\kern .6pt\box0}}}}

\def\morphism#1{\buildrel#1\over\longrightarrow}
\def\isomorphism#1{\mathrel{\mathop{\longrightarrow}%
\limits^{#1}_{\raise0.5ex\hbox{$\scriptstyle\sim$}}}}

\def\injlim{\mathop{\vtop{\offinterlineskip\halign{##\cr
 \hfil\rm lim\hfil\cr\noalign{\kern.1ex}\rightarrowfill\cr
 \noalign{\kern-.4ex}\cr}}}}
\def\projlim{\mathop{\vtop{\offinterlineskip\halign{##\cr
 \hfil\rm lim\hfil\cr\noalign{\kern.1ex}\leftarrowfill\cr
 \noalign{\kern-.4ex}\cr}}}}

\def\smallmatrix#1#2#3#4{\bigl({#1\atop#3}\,{#2\atop#4}\bigr)}
\def\medmatrix#1#2#3#4{\biggl({#1\atop#3}\,\,\,{#2\atop#4}\biggr)}

\def\blank{\mkern12mu}

\def\innerprod{\langle\blank\mathord,\blank\rangle}

\def\textfrac#1/#2{{\textstyle{#1\over#2}}}

\def\relativediag#1#2#3#4#5#6{\vcenter{\baselineskip=3ex \halign{
\hfil$##$&$##$&$##$\hfil\cr
#1\quad& \hfilneg\buildrel#2\over\longrightarrow\hfilneg& \quad#3\cr
\lower1ex\llap{$\scriptstyle#4\hskip-1ex$}\searrow& \quad&
\swarrow\lower1ex\rlap{$\hskip-1ex\scriptstyle#5$} \cr
& \hfilneg#6\hfilneg&\cr}}}

\def\trianglediag#1#2#3#4#5#6{\vcenter{\baselineskip=3ex \halign{
\hfil$##$&$##$&$##$\hfil\cr
#1\quad& \hfilneg\buildrel#2\over\longrightarrow\hfilneg& \quad#3\cr
\lower1ex\llap{$\scriptstyle#6\hskip-1ex$}\nwarrow& \quad&
\swarrow\lower1ex\rlap{$\hskip-1ex\scriptstyle#4$} \cr
& \hfilneg#5\hfilneg&\cr}}}

\def\correspondence#1#2#3#4#5{\vcenter{\baselineskip=3ex \halign{
\hfil$##$&$##$&$##$\hfil\cr
&\hfilneg#1\hfilneg\cr
\raise1ex\llap{$\scriptstyle#2$}\swarrow&&\searrow
\raise1ex\rlap{$\scriptstyle#3$}\cr
#4&&#5\cr}}}

\newif\iffirstpar
\everypar{\iffirstpar\parindent=\normalparindent\firstparfalse\fi}

\def\sectionheading#1{\subcount=0 \subsectioncount=0 \eqcount=0
  \bigskip\vskip\parskip
  \leftline{\bf #1}\nobreak\smallskip\firstpartrue\parindent=0pt}

\newif\ifappendix \appendixfalse

\def\currentsection{\ifappendix A\else\number\sectioncount\fi}

\def\section#1\par{\advance\sectioncount by1%
  \edef\currentlabel{\currentsection}%
  \sectionheading{\currentsection.\enspace#1}}

\def\unnumberedsection#1\par{\sectionheading{#1}}

\def\subsection#1\par{\medbreak\penalty-200\advance\subsectioncount by1%
  \edef\currentlabel{\currentsection.\number\subsectioncount}%
  \leftline{\it\currentsection.\number\subsectioncount.\enspace#1}%
  \smallskip\parindent=0pt\firstpartrue}

\newwrite\auxfile

\newcount\sectioncount \sectioncount=0
\newcount\subsectioncount 
\newcount\subcount 
\newcount\eqcount 

\def\subno{\global\advance\subcount by1\relax
  \currentsection.\number\subcount
  \xdef\currentlabel{\currentsection.\number\subcount}}
\def\proclaim #1. #2\par{\medbreak
  \noindent{\bf#1~\subno.\enspace}{\sl#2\par}%
  \ifdim\lastskip<\medskipamount \removelastskip\penalty55\medskip\fi}
\def\proclaimx #1 (#2). #3\par{\medbreak
  \noindent{\bf#1~\subno\ \rm (#2).\enspace}{\sl#3\par}%
  \ifdim\lastskip<\medskipamount \removelastskip\penalty55\medskip\fi}

\newdimen\algindent
\def\plusindent{\advance\algindent by \parindent}
\def\minusindent{\advance\algindent by-\parindent}


\newcount\algstepcount

\long\def\algorithm (#1). #2\endalgorithm{\medbreak
  \algindent=0pt%
  \algstepcount=0%
  \noindent{\bf Algorithm~\subno} (#1). {\sl#2}\par\medbreak}

\def\step{\advance\algstepcount by1
\edef\currentlabel{\number\algstepcount}
\smallskip\hangindent\parindent
\advance\hangindent by\algindent\indent
\llap{{\bf \the\algstepcount.}\enspace}\kern\algindent
\ignorespaces}


\def\labeldef#1#2{\expandafter\gdef\csname L@#1\endcsname{#2}}
\def\label#1{%
  \expandafter\xdef\csname L@#1\endcsname{\currentlabel}%
  \write\auxfile{\string\labeldef{#1}{\csname L@#1\endcsname}}%
  \ignorespaces}
\def\ref#1{{\rm\expandafter\ifx\csname L@#1\endcsname\relax
  \message{Undefined label `#1'}??\else
  \csname L@#1\endcsname\fi}}

\def\eqdef#1#2{\expandafter\gdef\csname E@#1\endcsname{#2}}
\def\eqnumber#1{\global\advance\eqcount by1\relax
  \eqno{\rm(\currentsection.\number\eqcount)}%
  \expandafter\xdef\csname E@#1\endcsname{%
    \currentsection.\number\eqcount}%
  \write\auxfile{\string\eqdef{#1}{\csname E@#1\endcsname}}}
\def\eqref#1{{\rm(\expandafter\ifx\csname E@#1\endcsname\relax
  \message{Undefined equation `#1'}??\else
  \csname E@#1\endcsname\fi)}}

\newcount\refcount \refcount=0
\def\citedef#1#2{\expandafter\gdef\csname C@#1\endcsname{#2}}
\def\cite#1{\expandafter\ifx\csname C@#1\endcsname\relax
  \message{Undefined reference `#1'}\citedef{#1}{??}\fi
  \expandafter\gdef\csname R@#1\endcsname{\relax}%
  [\csname C@#1\endcsname]}
\def\citex#1#2{\expandafter\ifx\csname C@#1\endcsname\relax
  \message{Undefined reference `#1'}\citedef{#1}{??}\fi
  \expandafter\gdef\csname R@#1\endcsname{\relax}%
  [\csname C@#1\endcsname, #2]}
\def\reference#1{\advance\refcount by 1%
  \expandafter\ifx\csname R@#1\endcsname\relax
  \message{Warning: reference `#1' not used}\fi
  \expandafter\edef\csname C@#1\endcsname{\the\refcount}%
  \write\auxfile{\string\citedef{#1}{\csname C@#1\endcsname}}%
  \item{[\csname C@#1\endcsname]}}

\newif\ifauxexists
\immediate\openin0=\jobname.aux
\ifeof 0
  \auxexistsfalse
\else
  \auxexiststrue
\fi
\immediate\closein0
\ifauxexists
  \input \jobname.aux
\else
  \message{No file `\jobname.aux'}
\fi
\openout\auxfile=\jobname.aux

\def\fc{{\frak c}}
\def\C{{\bf C}}
\def\fd{{\frak d}}
\def\genus_#1{g_{\lower1.6pt\hbox{$\scriptstyle#1$}}}
\def\gr{\mathop{\rm gr}\nolimits}
\def\grcan{\gr^{\rm can}}
\def\hyp{{\bf H}}
\def\Ltwo{{\rm L}^2}
\def\mucan{\mu^{\rm can}}
\def\muhyp{\mu_\hyp^{}}
\def\P{{\bf P}}
\def\R{{\bf R}}
\def\SL{{\rm SL}}
\def\vol{\mathop{\rm vol}\nolimits}
\def\X{{\rm X}}
\def\Z{{\bf Z}}

\font\titlefont=cmssbx10 scaled \magstep2

\centerline{\titlefont Explicit bounds on canonical Green functions of
modular curves}

\medskip
\centerline{Peter Bruin}
\smallskip
\centerline{25 July 2012}

\bigskip
{\narrower\noindent {\it Abstract\/.}\enspace 
We prove explicit bounds on canonical Green functions of Riemann
surfaces obtained as compactifications of quotients of the hyperbolic
plane by Fuchsian groups.\par}

\vfootnote{}{{\it Mathematics Subject Classification (2010):\/}
11F72,  
14G35,  
14G40,  
30F35,  
35J08   

This paper evolved from part of the author's thesis \cite{thesis}, the
research for which was supported by the Netherlands Organisation for
Scientific Research (NWO).  Further research was supported by Swiss
National Science Foundation grant 124737 and by the
Max-Planck-Institut f\"ur Mathematik, Bonn.}

\section Introduction

Let $\hyp$ denote the hyperbolic plane, identified with the complex
upper half-plane with holomorphic coordinate $z=x+iy$.  Let
$\muhyp=y^{-2}dx\,dy$ denote the standard volume form on~$\hyp$.  The
group $\SL_2(\R)$ acts isometrically on~$\hyp$ by $\gamma z={az+b\over
  cz+d}$ for $\gamma=\smallmatrix abcd\in\SL_2(\R)$ and $z\in\hyp$.
The Laplace operator on~$\hyp$ is
$\Delta=y^2(\partial_x^2+\partial_y^2)$.

A {\it Fuchsian group\/} is a discrete subgroup of~$\SL_2(\R)$.  A
Fuchsian group~$\Gamma$ is called {\it cofinite\/} if the volume of
$\Gamma\backslash\hyp$ with respect to the measure induced by~$\muhyp$
is finite.

Let $\Gamma$ be a cofinite Fuchsian group, and let $X$ be the standard
compactification of~$\Gamma\backslash\hyp$ obtained by adding the
cusps.  We assume that $X$ has positive genus.  There are two
interesting Green functions in this setting.  First, we have the
hyperbolic Green function~$\gr_\Gamma$ outside the diagonal
on~$\Gamma\backslash\hyp\times\Gamma\backslash\hyp$, given by the
structure of~$\Gamma\backslash\hyp$ as a quotient of the hyperbolic
plane by a Fuchsian group.  This appears naturally in the spectral
theory of automorphic forms.  Second, we have the canonical Green
function~$\grcan_X$ outside the diagonal on $X\times X$, given by the
structure of~$X$ as a compact Riemann surface of positive genus.  This
is a fundamental object in the intersection theory on arithmetic
surfaces developed by Arakelov~\cite{Arakelov},
Faltings~\cite{Faltings: Calculus on arithmetic surfaces} and others,
where it is used to define local intersection numbers of horizontal
divisors at the infinite places.


In this article we derive explicit bounds on~$\grcan_X$ using bounds
on~$\gr_\Gamma$ established by the author in~\cite{green1}.  Our
results are valid for any cofinite Fuchsian group, although we are
motivated by the case of arithmetic groups, and in particular
congruence subgroups of~$\SL_2(\Z)$.  Bounds on the canonical Green
functions of the modular curves $\X_1(n)$ are relevant to recent work
of Edixhoven, Couveignes et~al.~\cite{book} and of the author
\cite{thesis}, where Arakelov theory is employed to obtain a
polynomial-time algorithm for computing Galois representations
attached to Hecke eigenforms over finite fields.

The following theorem illustrates our general results.  To avoid
having to deal with the logarithmic singularity, we only give an upper
bound.  More precise results are contained in
Theorem~\ref{theorem:bounds-grcan}.

\proclaim Theorem. Let $\Gamma$ be a congruence subgroup of level~$n$
of\/~$\SL_2(\Z)$ such that the compactification~$X$ of\/
$\Gamma\backslash\hyp$ has positive genus.  Then the canonical Green
function of~$X$ satisfies
$$
\sup_{X\times X}\grcan_X\le 1.6\cdot 10^4+7.7n+0.088n^2.
$$

\label{theorem:example}

It will make little difference to us whether the quotient
$\Gamma\backslash\hyp$ and its compactification~$X$ are interpreted as
as stacks with generic stabiliser $\Gamma\cap\{\pm1\}$ or as their
coarse moduli spaces, which are Riemann surfaces.  To avoid any
subtleties, the reader can restrict himself to groups containing
neither $-1$ nor any elliptic elements and with all cusps regular,
such as $\Gamma_1(n)$ with $n\ge5$.  Two points are useful to keep in
mind.  First, as in~\cite{green1}, we define integration on
$\Gamma\backslash\hyp$ and~$X$ in a ``stack-like'' way, so that in
case $-1\in\Gamma$, integrals over $\Gamma\backslash\hyp$ or~$X$ are
half of what the na{\"\i}ve definition gives.  Second, the space of
cusp forms of weight~2 for~$\Gamma$, the space of holomorphic
differentials on~$X$, and the space of holomorphic differentials on
the coarse moduli space of~$X$ are all isomorphic.  We write
$\genus_X$ for the dimension of these spaces, and call it the genus
of~$X$.

\remark A different approach to the problem of bounding canonical
Green functions was taken by Jorgenson and
Kramer~\cite{Jorgenson-Kramer: Green's functions}.  For {\it
compact\/} quotients of the upper half-plane, they deduced an
interesting expression for the canonical Green function in terms of
data associated to the hyperbolic metric.  In comparison, our method
is less involved and applies more naturally to modular curves.

\goodbreak

\section Notation and statement of results

\subsection Cusps

\label{subsec:cusps}

Let $\Gamma$ be a cofinite Fuchsian group.  The cusps of~$\Gamma$
correspond to the conjugacy classes of maximal parabolic subgroups
of~$\Gamma$.  For every cusp $\fc$ of~$\Gamma$, we fix one such
subgroup and denote it by~$\Gamma_\fc$, and we fix
$\sigma_\fc\in\SL_2(\R)$ such that
$$
\{\pm1\}
\sigma_\fc^{-1}\Gamma_\fc\sigma_\fc=
\biggl\{\pm\medmatrix 1b01\biggm|b\in\Z\biggr\}.
$$
For $z\in\hyp$, we write
$$
q_\fc(z)=\exp(2\pi i\sigma_\fc^{-1}z)
$$
and
$$
y_\fc(z)=\Im\sigma_\fc^{-1}z
=-{\log|q_\fc(z)|\over2\pi}.
$$

For all $\gamma\in\Gamma$, we write
$$
C_\fc(\gamma)=|c|\quad\hbox{if }
\sigma_\fc^{-1}\gamma\sigma_\fc=\medmatrix abcd.
$$
It is known that if $\epsilon$ is a real number satisfying
the inequality
$$
0<\epsilon\le\min_{\textstyle{\gamma\in\Gamma\atop
\gamma\not\in\Gamma_\fc}}C_\fc(\gamma),
\eqnumber{ineq:eps-sufficiently-small}
$$
then for all $z\in\hyp$ and~$\gamma\in\Gamma$ one has the implication
$$
y_\fc(z)>1/\epsilon\hbox{ and }y_\fc(\gamma z)>1/\epsilon
\;\Longrightarrow\;\gamma\in\Gamma_\fc.
$$
For any $\epsilon$ satisfying \eqref{ineq:eps-sufficiently-small}, the
image of the strip
$$
\{x+iy\mid 0\le x<1\hbox{ and }y>1/\epsilon\}\subset\hyp
$$
under the map
$$
\hyp\morphism{\sigma_\fc}\hyp\longrightarrow\Gamma\backslash\hyp
$$
is an open disc $D_\fc(\epsilon)$ around~$\fc$.  The map~$q_\fc$
induces a chart on~$\Gamma\backslash\hyp$ identifying
$D_\fc(\epsilon)$ with the punctured disc $\{z\in\C\mid
0<|z|<\exp(-2\pi/\epsilon)\}$.  The disc~$D_\fc(\epsilon)$ has the
compactification
$$
\bar D_\fc(\epsilon)=\{z\in\Gamma\backslash\hyp\mid
y_\fc(z)\ge 1/\epsilon\}\cup\{\fc\}.
$$

\subsection The canonical $(1,1)$-form

Let $X$ be a compact connected Riemann surface of genus
$\genus_X\ge1$.  The $\C$-vector space $\Omega^1(X)$ of global
holomorphic differentials on~$X$ has dimension~$\genus_X$ and is
equipped with the inner product
$$
\langle\alpha,\beta\rangle={i\over2}\int_X\alpha\wedge\bar\beta.
$$
The {\it canonical $(1,1)$-form on~$X$\/} is
$$
\mucan_X={i\over2\genus_X}\sum_{\alpha\in B}\alpha\wedge\bar\alpha,
$$
where $B$ is any orthonormal basis of~$\Omega^1(X)$ with respect
to~$\innerprod$.  The form~$\mucan_X$ is independent of the choice
of~$B$.

Let us now assume that $X$ is (the coarse moduli space associated to)
the compactification of $\Gamma\backslash\hyp$ with $\Gamma$ a
cofinite Fuchsian group.  We define a smooth and bounded
function~$F_{\Gamma}$ on~$\Gamma\backslash\hyp$ by
$$
F_\Gamma(z)=\sum_{f\in B}(\Im z)^2|f(z)|^2,
$$
where $B$ is any orthonormal basis for the space of holomorphic cusp
forms of weight~2 for~$\Gamma$.  The $(1,1)$-forms $\mucan_X$
and~$\muhyp$ are related by
$$
\mucan_X={1\over\genus_X}F_\Gamma\muhyp.
\eqnumber{eq:mucan-FGamma}
$$

\subsection Spectral theory for Fuchsian groups

We collect here some facts that we will need.  For proofs and further
details, we refer to Iwaniec's book \cite{Iwaniec} or Hejhal's two
volumes \cite{Hejhal-I} and~\cite{Hejhal-II}.

Let $\Gamma$ be a cofinite Fuchsian group.  The Laplace operator
on~$\hyp$ induces an (unbounded, densely defined) self-adjoint
operator $\Delta_\Gamma$ on the Hilbert space
$\Ltwo(\Gamma\backslash\hyp)$.  We will often work with
$-\Delta_\Gamma$ instead, since it is non-negative.  The spectrum
of~$-\Delta_\Gamma$ consists of a discrete part and a continuous part.
The discrete part consists of eigenvalues, which we we denote by
$$
0=\lambda_0<\lambda_1\le\lambda_2\le\cdots,
\quad\lambda_j\to\infty\hbox{ as }j\to\infty.
$$
Let $\{\phi_j\}_{j=0}^\infty$ be a corresponding orthonormal system of
eigenfunctions.  The continuous part of the spectrum
of~$-\Delta_\Gamma$ is the interval $[1/4,\infty)$ with multiplicity
equal to the number of cusps of~$\Gamma$.  It corresponds to the
non-holomorphic Eisenstein series $E_\fc(z,s)$.

Every smooth and bounded function on $\Gamma\backslash\hyp$ has the
spectral representation
$$
f(z)=\sum_{j=0}^\infty b_j\phi_j(z)
+\sum_\fc{1\over 4\pi i}\int_{\Re s=1/2} b_\fc(s)\,
E_\fc(z,s)ds,
$$
where $\fc$ runs over the cusps of~$\Gamma$ and the coefficients $b_j$
and $b_\fc(s)$ are given by
$$
b_j=\int_{\Gamma\backslash\hyp}f\bar\phi_j\muhyp
\quad\hbox{and}\quad
b_\fc(s)=\int_{\Gamma\backslash\hyp}f
\bar E_\fc(\blank, s)\muhyp.
$$

For later use, we define, as in \citex{green1}{\S2.5},
$$
\Phi_\Gamma(z,\lambda)=
\sum_{j\colon\,\lambda_j\le\lambda}|\phi_j(z)|^2
+\sum_\fc{1\over 4\pi i}
\int_{\textstyle{\Re s=1/2\atop s(1-s)\le\lambda}}
\bigl|E_\fc(z,s)\bigr|^2ds.
\eqnumber{eq:def-PhiGamma}
$$
It is known that for $\lambda$ tending to~$\infty$ and fixed $z$, this
function is bounded linearly in~$\lambda$.

\subsection The Green function of a Fuchsian group

Let $\Gamma$ be a cofinite Fuchsian group.  For every smooth and
bounded function~$f$ on~$\Gamma\backslash\hyp$, there exists a unique
smooth and bounded function $g_f$ on~$\Gamma\backslash\hyp$ such that
$$
\Delta_\Gamma g_f=f-{1\over \vol_\Gamma}\int_{\Gamma\backslash\hyp}f\muhyp
\quad\hbox{and}\quad\int_{\Gamma\backslash\hyp}g_f\muhyp=0.
$$
There exists a unique function $\gr_\Gamma$ on
$\Gamma\backslash\hyp\times\Gamma\backslash\hyp$ that satisfies
$\gr_\Gamma(z,w)=\gr_\Gamma(w,z)$, is smooth except for a logarithmic
singularity along the diagonal, and has the property that if $f$ is a
smooth and bounded function on~$\Gamma\backslash\hyp$, then the
function~$g_f$ is given by
$$
g_f(z)=\int_{w\in\Gamma\backslash\hyp}\gr_\Gamma(z,w)f(w)\muhyp(w).
$$
The function~$\gr_\Gamma$ is called the {\it Green function of the
Fuchsian group~$\Gamma$\/}.

\subsection The canonical Green function of a Riemann surface

Let $X$ be a Riemann surface.  Let $*$ denote the star operator on
smooth $1$-forms, given with respect to any local holomorphic
coordinate~$z=x+iy$ by
$$
*dx=dy,\quad *dy=-dx.
$$
If we identify $X$ locally with the hyperbolic plane, an easy
calculation shows that the operator $d*d$ sending functions to
$(1,1)$-forms is related to the Laplace operator~$\Delta$ as follows:
if $f$ is any smooth function on~$X$, then
$$
d*df=\Delta f\cdot\muhyp.
$$

Let us now assume that $X$ is compact, connected and of positive
genus.  Let $\alpha$ be a smooth $(1,1)$-form on~$X$.  Then there
exists a unique smooth function $h_f$ on~$X$ such that
$$
d*d h_\alpha=\alpha-\biggl(\int_X\alpha\biggr)\mucan_X
\quad\hbox{and}\quad\int_X h_\alpha\mucan_X=0.
$$
There exists a unique function $\grcan_X$ on $X\times X$ that
satisfies $\grcan_X(z,w)=\grcan_X(w,z)$, is smooth except for a
logarithmic singularity along the diagonal, and has the property that
if $\alpha$ is a smooth $(1,1)$-form on~$X$, then the
function~$h_\alpha$ is given by
$$
h_\alpha(z)=\int_{w\in X}\grcan_X(z,w)\alpha(w).
$$
The function~$\grcan_X$ is called the {\it canonical Green function\/}
of~$X$.

\subsection The main result

For all $u>1$, we write
$$
L(u)={1\over4\pi}\log{u+1\over u-1}.
$$
If $\Gamma$ is a Fuchsian group and $X$ is the compactification
of~$\Gamma\backslash\hyp$ obtained by adding the cusps, we write
$$
\zeta_\Gamma={1\over\genus_X}
\int_{\Gamma\backslash\hyp}F_\Gamma\mucan_X-{1\over\vol_\Gamma}.
\eqnumber{eq:definition-zeta}
$$
We will commit the following abuse of notation: for $z\in\hyp$ and $Z$
a subset of~$\Gamma\backslash\hyp$, we write $z\in Z$ if the image
of~$z$ in~$\Gamma\backslash\hyp$ lies in~$Z$.

\proclaim Theorem.  Let $\Gamma$ be a cofinite Fuchsian group, and let
$X$ be the compactification of\/ $\Gamma\backslash\hyp$ obtained by
adding the cusps.  Let $\delta$ be a real number with $\delta>1$.  For
every cusp~$\fc$ of\/~$\Gamma$, let $\epsilon_\fc'>\epsilon_\fc>0$ be
real numbers satisfying the inequalities
$$
\epsilon_\fc'\bigl(\delta+\sqrt{\delta^2-1}\bigr)^{1/2}\le
\min_{\textstyle{\gamma\in\Gamma\atop\gamma\not\in\Gamma_\fc}}C_\fc(\gamma)
\quad\hbox{and}\quad
\bigl(\delta+\sqrt{\delta^2-1}\bigr)\epsilon_\fc\le\epsilon_\fc'
$$
and small enough such that the discs $D_\fc(\epsilon_\fc')$ are
pairwise disjoint.  Let $Y$ be the compact subset
of~$\Gamma\backslash\hyp$ defined by
$$
Y=(\Gamma\backslash\hyp)\mathbin{\big\backslash}
\bigsqcup_{\fc\;\rm cusp}D_\fc(\epsilon_\fc).
$$
Let $A$ and~$B$ be real numbers such that the hyperbolic Green
function~$\gr_\Gamma$ satisfies
$$
A\le
\gr_\Gamma(z,w)+\sum_{\textstyle{
\gamma\in\Gamma\atop u(z,\gamma w)\le\delta}}
\bigl(L(u(z,\gamma w)-L(\delta))\bigr)
\le B
\quad\hbox{for all }z,w\in Y.
\eqnumber{ineq:grGammaAB}
$$
Let $C>0$ be such that the function~$\Phi_\Gamma(z,\lambda)$
defined by~\eqref{eq:def-PhiGamma} satisfies
$$
\Phi_\Gamma(z,\lambda)\le C\lambda
\quad\hbox{for all $z\in Y$ and $\lambda\ge1/4$}.
\eqnumber{ineq:C}
$$
Let $\eta\in(0,1/4]$ be such that the spectrum of~$-\Delta_\Gamma$ is
contained in $\{0\}\cup[\eta,\infty)$.  With the notation
$$
\eqalign{
S&=\sqrt{\biggl({1\over4\eta^2}+4\biggr)\,C\zeta_\Gamma},\cr
T(\epsilon)&={\sup_Y F_\Gamma\over\genus_X}
\Bigl({\epsilon\over4\pi}\Bigr)^2\quad\hbox{for all }\epsilon>0,\cr
r_\delta&={1\over24\pi}\biggl(\sqrt{2\over\delta-1}
+\arctan\sqrt{\delta-1\over2}\biggr),\cr
\tilde A_\fc&=A
+\#(\Gamma\cap\{\pm1\})\biggl[{1\over\epsilon_\fc'}
\biggl(1-{2\over\pi}\arctan\sqrt{\delta-1\over2}\biggr)
-\epsilon_\fc' r_\delta\biggr],\cr
\tilde B_\fc&=B
+\#(\Gamma\cap\{\pm1\})\biggl[{1\over\epsilon_\fc'}
\biggl(1-{2\over\pi}\arctan\sqrt{\delta-1\over2}\biggr)
+\epsilon_\fc' r_\delta\biggr],}
$$
we have the following bounds on the canonical Green function
$\grcan_X(z,w)$:
\smallskip\noindent
(a)\enspace If $z,w\in Y$, we have
$$
A-2S-\zeta_\Gamma/\eta\le
\grcan_X(z,w)
+\sum_{\textstyle{\gamma\in\gamma\atop u(z,\gamma w)\le\delta}}
\bigl(L(u(z,w)-L(\delta))\bigr)
\le B+2S.
$$
(b)\enspace If $\fc$ is a cusp such that $z\in D_\fc(\epsilon_\fc)$,
$w\in Y$ and $w\not\in D_\fc(\epsilon_\fc')$, or such that $w\in
D_\fc(\epsilon_\fc)$, $z\in Y$ and $z\not\in D_\fc(\epsilon_\fc')$,
then we have
$$
A-2S-\zeta_\Gamma/\eta\le\grcan_X(z,w)\le B+2S
+T(\epsilon_\fc).
$$
(c)\enspace If $\fc$, $\fd$ are two distinct cusps such that $z\in
D_\fc(\epsilon_\fc)$ and $w\in D_\fd(\epsilon_\fd)$, we have
$$
A-2S-\zeta_\Gamma/\eta\le\grcan_X(z,w)\le B+2S
+T(\epsilon_\fc)+T(\epsilon_\fd).
$$
(d)\enspace If $\fc$ is a cusp such that $z,w\in
D_\fc(\epsilon_\fc')$, we have
$$
\tilde A_\fc-2S-\zeta_\Gamma/\eta\le\grcan_X(z,w)
-\#(\Gamma\cap\{\pm1\})\cdot{1\over2\pi}\log|q_\fc(z)-q_\fc(w)|
\le\tilde B_\fc+2S+2T(\epsilon_\fc').
$$

\label{theorem:bounds-grcan}

In earlier work of the author~\cite{green1}, it was described how to
compute explicit real numbers $A$ and~$B$ as in~\eqref{ineq:grGammaAB}
and $C$ as in~\eqref{ineq:C} for concrete groups~$\Gamma$, such as
congruence groups of~$\SL_2(\Z)$.
In Section~\ref{sec:bounds-cusp-forms} below, we will show that bounds
on the function~$F_\Gamma$ can likewise be found easily for given
groups~$\Gamma$.

\section Tools

\subsection Comparison of hyperbolic and canonical Green functions

There is a standard way to relate the hyperbolic and canonical Green
functions, which we will use to find explicit bounds on the canonical
Green function.  Let $\Gamma$ be a cofinite Fuchsian group, and let
$X$ be the compactification of~$\Gamma\backslash\hyp$.  We define a
function~$h_\Gamma\colon\Gamma\backslash\hyp\to\R$ by
$$
\eqalign{
h_\Gamma(z)&=\int_{w\in\Gamma\backslash\hyp}\gr_\Gamma(z,w)\mucan_X(w)\cr
&={1\over\genus_X}\int_{w\in\Gamma\backslash\hyp}
\gr_\Gamma(z,w)F_\Gamma(w)\muhyp(w).}
\eqnumber{eq:hGamma}
$$
By the defining properties of~$\gr_\Gamma$, the function~$h_\Gamma$
satisfies
$$
\eqalign{
\Delta h_\Gamma&={1\over\genus_X}F_\Gamma
-{1\over\genus_X\vol_\Gamma}\int_{\Gamma\backslash\hyp}F_\Gamma\muhyp\cr
&={1\over\genus_X}F_\Gamma-{1\over\vol_\Gamma}.}
$$
This implies that the canonical Green function of~$X$ can be expressed
as
$$
\grcan_X(z,w)=\gr_\Gamma(z,w)-h_\Gamma(z)-h_\Gamma(w)
+\int_{\Gamma\backslash\hyp} h_\Gamma\mucan_X.
\eqnumber{comparison-green}
$$

\subsection The Selberg--Harish-Chandra transform

\label{subsec:shc}

Let $k$ be a real number, and let
$\Delta_k=y^2(\partial_x^2+\partial_y^2)-iky\partial_x$ denote the
Laplace operator of weight~$k$ on~$\hyp$.  Let
$\theta\colon[1,\infty)\to\R$ be a piecewise smooth function with
compact support.  We define
$$
\theta^{(k)}(z,w)=
\biggl({w-\bar z\over z-\bar w}\biggr)^{k/2}\theta(u(z,w)).
$$

Let $P_{s,k}$ be the generalisation of the Legendre
function~$P_{s-1}(u)$ given by Fay \citex{Fay}{\S1} (note that our
definition of weight is twice that of~\cite{Fay}):
$$
P_{s,k}=\biggl({2\over u+1}\biggr)^s
F\biggl(s-{k\over2},s+{k\over2};1;{u-1\over u+1}\biggr).
\eqnumber{eq:Psk}
$$
We define the {\it Selberg--Harish-Chandra transform\/} of weight~$k$
of the function~$\theta$ as
$$
h_\theta^{(k)}(s)=2\pi\int_1^\infty\theta(u)P_{s,k}(u)du.
\eqnumber{eq:h-formula}
$$
If $f$ is an eigenfunction of~$-\Delta_k$ with eigenvalue $s(1-s)$,
then we have
$$
\int_{w\in\hyp}\theta^{(k)}(z,w) f(w) dw=h_\theta^{(k)}(s)f(z);
\eqnumber{eq:shc}
$$
see Fay \citex{Fay}{Theorem 1.5}.

\subsection Automorphic forms

\label{subsec:automorphic-forms}

For simplicity, we take $k=2$ from now on.  We write
$$
\nu(\gamma,z)={cz+d\over c\bar z+d}
={(cz+d)^2\over|cz+d|^2}
\quad\hbox{for }\gamma=\medmatrix abcd\in\SL_2(\R)
\hbox{ and }z\in\hyp.
$$

We recall that an {\it automorphic form (of Maa\ss)\/} of weight~2
for~$\Gamma$ is a smooth function~$f\colon\hyp\to\C$ with the
following properties:
\smallskip
\item{(1)} the function~$f$ satisfies the transformation formula
$$
f(\gamma z)=\nu(\gamma,z)f(z)
\quad\hbox{ for all $\gamma\in\Gamma$ and $z\in\hyp$};
$$
\smallskip
\item{(2)} for every cusp~$\fc$ of~$\Gamma$, there is a real
number~$\kappa$ such that $|f(z)|=O(y_\fc(z)^\kappa)$ as
$y_\fc(z)\to\infty$.
\smallskip\noindent
A {\it cusp form\/} of weight~2 for~$\Gamma$ is a function~$f$
satisfying (1) and the following condition (which is stronger
than~(2)):
\smallskip
\item{($2'$)} for every cusp~$\fc$ of~$\Gamma$ there exists
$\epsilon>0$ such that $|f(z)|=O(\exp(-\epsilon y_\fc z))$ as
$y_\fc(z)\to\infty$.

\smallskip
Let $\Ltwo(\Gamma\backslash\hyp,2)$ denote the Hilbert space of
square-integrable automorphic forms of weight~2 for~$\Gamma$, equipped
with the Petersson inner product.

Let $\theta$ be a function as in~\S\ref{subsec:shc}.  Then we have
$$
\theta^{(2)}(\gamma z,\gamma w)={\nu(\gamma,z)\over\nu(\gamma,w)}
\theta^{(2)}(z,w)
\quad\hbox{for all $\gamma\in\SL_2(\R)$ and $z,w\in\hyp$}.
$$
Let $\Gamma$ be a cofinite Fuchsian group.  We define
$$
K_{\Gamma,\theta}^{(2)}(z,w)=\sum_{\gamma\in\Gamma}
\nu(\gamma,w)\theta^{(2)}(z,\gamma w).
\eqnumber{eq:def-K}
$$
This function satisfies
$$
K_{\Gamma,\theta}^{(2)}(w,z)=\overline{K_{\Gamma,\theta}^{(2)}(z,w)}
$$
and, for all $\gamma\in\Gamma$,
$$
\eqalign{
K_{\Gamma,\theta}^{(2)}(\gamma z,w)&=\nu(\gamma,z)K_{\Gamma,\theta}^{(2)}(z,w),\cr
K_{\Gamma,\theta}^{(2)}(z,\gamma w)&=\nu(\gamma,w)^{-1}K_{\Gamma,\theta}^{(2)}(z,w).}
$$
Now \eqref{eq:shc} implies that if $f$ is an automorphic form of
weight~2 for~$\Gamma$ satisfying $-\Delta_2 f=s(1-s)f$, then
$$
\int_{w\in\Gamma\backslash\hyp}K_{\Gamma,\theta}^{(2)}(z,w) f(w) \muhyp(w)
=h_\theta^{(2)}(s)f(z).
\eqnumber{eq:shc-Gamma}
$$

\section Explicit bounds on the canonical $(1,1)$-form

\label{sec:bounds-cusp-forms}

Let $\Gamma$ be a cofinite Fuchsian group.  In this section we find
bounds on the function~$F_\Gamma$ that are easy to evaluate explicitly
in concrete cases.  We essentially adapt the methods of
Iwaniec~\citex{Iwaniec}{\S7.2} from weight~0 to weight~2.  This method
is more elementary than that used by Jorgenson and Kramer in
\cite{Jorgenson-Kramer: automorphic forms},
and our bounds are easy to make explicit, as the example in
Section~\ref{sec:example} shows.

For $z\in\hyp$ and $b\ge1$, we write
$$
N_\Gamma(z,b)=\#\{\gamma\in\Gamma\mid u(z,\gamma z)\le b\}.
$$

\proclaim Proposition.  For every cofinite Fuchsian group $\Gamma$,
all $z\in\hyp$ and all $a>1$, we have
$$
F_\Gamma(z)\le{(a-1)N_\Gamma(z,2a^2-1)\over
8\pi\bigl(\log{a+1\over2}\bigr)^2}.
$$

\label{prop:bound-F}

\proof Let $(f_1,\ldots,f_g)$ be an orthonormal basis of the space of
holomorphic cusp forms of weight~2 for~$\Gamma$.  We write
$$
\phi_j(z)=(\Im z)f_j(z).
$$
Then the $\phi_j$ are annihilated by the operator~$\Delta_2$, and
$(\phi_1,\ldots,\phi_g)$ is an orthonormal system in the Hilbert space
$\Ltwo(\Gamma\backslash\hyp,2)$ of automorphic forms of weight~2
for~$\Gamma$.

Let $z\in\hyp$ and $a>1$ be given.  We apply \SS\ref{subsec:shc}
and~\ref{subsec:automorphic-forms} with
$$
\theta(u)=\cases{
1& if $1\le u\le a$,\cr
0& if $u>a$.}
$$
In the Hilbert space $\Ltwo(\Gamma\backslash\hyp,2)$, we consider
$\overline{K_{\Gamma,\theta}^{(2)}(z,w)}$, as a function of~$w$, and
the orthonormal system $(\phi_1,\ldots,\phi_g)$.  From Bessel's
inequality and~\eqref{eq:shc-Gamma}, we obtain
$$
\sum_{j=1}^g\bigl|h_\theta^{(2)}(0)\phi_j(z)\bigr|^2\le
\int_{w\in\Gamma\backslash\hyp}\bigl|K_{\Gamma,\theta}^{(2)}(z,w)\bigr|^2\muhyp(w).
$$
We note that the left-hand side is equal to
$\bigl|h_\theta^{(2)}(0)\bigr|^2F_\Gamma(z)$.  Let us denote the
right-hand side by $\kappa(z)$; this is a $\Gamma$-invariant function
of~$z$ with values in $[0,\infty)$.  The definition~\eqref{eq:def-K}
gives
$$
\kappa(z)=\sum_{\gamma_1,\gamma_2\in\Gamma}
\nu(\gamma_1,w)\theta^{(2)}(z,\gamma_1 w)
\overline{\nu(\gamma_2,w)\theta^{(2)}(z,\gamma_2 w)}\muhyp(w).
$$
Putting $\gamma=\gamma_1\gamma_2^{-1}$, we obtain after a
straightforward computation
$$
\kappa(z)=\sum_{\gamma\in\Gamma}\nu(\gamma,z)\int_{w\in\hyp}
\theta^{(2)}(z,w)\overline{\theta^{(2)}(\gamma z,w)}\muhyp(w).
$$
This implies
$$
\kappa(z)\le\sum_{\gamma\in\Gamma}\int_{w\in\hyp}
\theta(u(z,w))\theta(u(\gamma z,w))\muhyp(w).
$$
By the definition of~$\theta$, the integral on the right-hand side can
be interpreted as the area of the intersection of the discs of area
$2\pi(a-1)$ around $z$ and~$\gamma z$, respectively.  By the triangle
area for the hyperbolic distance, this intersection is empty unless
$$
u(z,\gamma z)\le 2a^2-1.
$$
This implies
$$
\kappa(z)\le2\pi(a-1)N_\Gamma(z,2a^2-1),
$$
and hence
$$
|h_\theta^{(2)}(0)|^2 F_\Gamma(z)\le 2\pi(a-1)N_\Gamma(z,2a^2-1).
$$

We evaluate $h_\theta^{(2)}(0)$ using \eqref{eq:Psk}
and~\eqref{eq:h-formula}.  The hypergeometric series terminates after
two terms and gives
$$
P_{0,2}(u)={2\over u+1}.
$$
This implies
$$
h_\theta^{(2)}(0)=4\pi\log{u+1\over 2}.
$$
This finishes the proof.\endproof

The above proposition does not give the correct asymptotic behaviour
of~$F_\Gamma(z)$ for $z$ close to a cusp of~$\Gamma$.  The following
result extends our bounds to neighbourhoods of the cusps.

\proclaimx Lemma (cf.\ Jorgenson and Kramer \citex{Jorgenson-Kramer:
  automorphic forms}{Theorem~3.1}). Let $\Gamma$ be a cofinite
Fuchsian group, let $\fc$ be a cusp of\/~$\Gamma$, and let $\epsilon$
be a real number satisfying~\eqref{ineq:eps-sufficiently-small}.  Then
for all $z\in D_\fc(\epsilon)$, we have
$$
\eqalign{
F_\Gamma(z)&\le(\epsilon y_\fc(z))^2
\exp(4\pi/\epsilon-4\pi y_\fc(z))
\sup_{\partial\bar D_\fc(\epsilon)} F_\Gamma\cr
&\le\cases{\displaystyle
\sup_{\partial\bar D_\fc(\epsilon)} F_\Gamma
& if $\epsilon\le2\pi$,\cr
\displaystyle
\Bigl({\epsilon\over 2\pi}\exp(2\pi/\epsilon-1)\Bigr)^2
\sup_{\partial\bar D_\fc(\epsilon)} F_\Gamma
& if $\epsilon>2\pi$.}}
$$

\label{lemma:extension-cusps-F}

\proof
Every holomorphic cusp form $f$ of weight~2 for~$\Gamma$ has a
$q$-expansion of the form
$$
f(z)dz=\sum_{n=1}^\infty a_{\fc,n}(f)q_\fc(z)^n
\cdot d(\sigma_\fc^{-1}z)
\quad\hbox{with }a_{\fc,n}(f)\in\C.
$$
This implies
$$
(\Im z)^2|f(z)|^2=y_\fc(z)^2\biggl|\sum_{n=1}^\infty
a_{\fc,n}(f)q_\fc(z)^n\biggr|^2.
$$
Applying this to an orthonormal basis of the space of holomorphic cusp
forms of weight~2, we see that the function
$$
y_\fc(z)^{-2}\exp(4\pi y_\fc(z))F_\Gamma(z)
=\sum_{f\in B}\left|{f(z)\over q_\fc(z)}\right|^2
$$
extends to a subharmonic function on $\bar D_\fc(\epsilon)$.  By the
maximum principle for subharmonic functions, the function assumes its
maximum on the boundary.  This implies the first inequality.  The
second inequality follows from the easily checked fact that the
function $(\epsilon y_\fc(z))^2\exp(4\pi/\epsilon-4\pi y_\fc(z))$ for
$y\ge1/\epsilon$ assumes its maximum at $y=1/(2\pi)$ if $\epsilon>
2\pi$, and at $y=1/\epsilon$ if $\epsilon\le 2\pi$.\endproof

\section Proof of the main result

Our proof of Theorem~\ref{theorem:bounds-grcan} is based on the
equation~\eqref{comparison-green}, the bounds on the hyperbolic Green
function from~\cite{green1}, and on bounds on the function~$h_\Gamma$
defined by~\eqref{eq:hGamma}.  The proof of the latter bounds occupies
most of this section; Theorem~\ref{theorem:bounds-grcan} then follows
without difficulties.

\proclaim Lemma. Let $\Gamma$, $Y$, $\eta$ and $C$ be as in
Theorem~\ref{theorem:bounds-grcan}.  Then the function
$$
M_\Gamma(z)=\sum_{j\ge1}{1\over\lambda_j^2}|\phi_j(z)|^2
+\sum_\fc{1\over 4\pi i}\int_{\Re s=1/2}{1\over(s(1-s))^2}
|E_\fc(z,s)|^2ds
$$
satisfies
$$
M_\Gamma(z)\le\biggl({1\over4\eta^2}+4\biggr)C
\quad\hbox{for all }z\in Y.
$$

\label{lemma:MGamma}

\proof Separating the terms with $\lambda_j\le1/4$, we get (with
$\partial\Phi_\Gamma/\partial\lambda$ taken in a distributional sense)
$$
\eqalign{
M_\Gamma(z)&=\sum_{j\colon\,0<\lambda_j\le1/4}
{1\over\lambda_j^2}|\phi_j(z)|^2
+\int_{1/4}^\infty{1\over\lambda^2}
{\partial\Phi_\Gamma(z,\lambda)\over\partial\lambda}d\lambda\cr
&\le{1\over\eta^2}\Phi(z,1/4)+\Bigl[{1\over\lambda^2}
\Phi_{\Gamma}(z,\lambda)\Bigr]_{\lambda=1/4}^\infty
+2\int_{1/4}^\infty\lambda^{-3}\Phi_\Gamma(z,\lambda)d\lambda\cr
&=\biggl({1\over\eta^2}-16\biggr)\Phi_\Gamma(z,1/4)
+2\int_{1/4}^\infty\lambda^{-3}\Phi_\Gamma(z,\lambda)d\lambda\cr
&\le\biggl({1\over\eta^2}-16\biggr){C\over4}
+2C\int_{1/4}^\infty\lambda^{-2}d\lambda\cr
&=\biggl({1\over\eta^2}-16\biggr){C\over4}+8C\cr
&=\biggl({1\over4\eta^2}+4\biggr)C,}
$$
where the second inequality follows from~\eqref{ineq:C}.\endproof

\proclaim Lemma. Let $\Gamma$, $Y$, $\eta$ and~$C$ be as in
Theorem~\ref{theorem:bounds-grcan}, and let $\zeta_\Gamma$ be as
in~\eqref{eq:definition-zeta}.  Then we have
$$
|h_\Gamma(z)|^2\le\biggl({1\over4\eta^2}+4\biggr)C\zeta_\Gamma
\quad\hbox{for all }z\in Y.
$$

\label{lemma:bound-hGamma-Y}

\proof Let $X$ be the compactification of~$\Gamma\backslash\hyp$.
Since the function~$F_\Gamma$ is smooth and bounded, we may consider
its spectral representation, say
$$
{1\over\genus_X}F_\Gamma(z)=\sum_{j\ge0} b_j\phi_j(z)
+\sum_\fc{1\over 4\pi i}\int_{\Re s=1/2}b_\fc(s)E_\fc(z,s)ds.
\eqnumber{eq:FGamma-spectral}
$$
Then the definition of~$h_\Gamma$ implies that it has the spectral
representation
$$
h_\Gamma(z)=-\sum_{j\ge 1}{b_j\over\lambda_j}\phi_j(z)
-\sum_\fc{1\over 4\pi i}\int_{\Re s=1/2}
{b_\fc(s)\over s(1-s)}E_\fc(z,s)ds.
\eqnumber{eq:hGamma-spectral}
$$
(Note the absence of the term corresponding to $j=0$.)  Now the
Cauchy--Schwarz inequality implies
$$
h_\Gamma(z)^2\le M_\Gamma(z)\Biggl(\sum_{j\ge1}|b_j|^2
+\sum_\fc{1\over 4\pi i}\int_{\Re s=1/2}|b_\fc(s)|^2\Biggr).
$$
Next, it follows from \eqref{eq:FGamma-spectral}, the identity
$|a_0^2|=1/{\vol_\Gamma}$ and~\eqref{eq:mucan-FGamma} that
$$
\eqalign{
\sum_{j\ge1}|b_j|^2
+\sum_\fc{1\over 4\pi i}\int_{\Re s=1/2}|b_\fc(s)|^2
&=\int_{z\in\Gamma\backslash\hyp}\biggl({1\over\genus_X}F_\Gamma(z)
-{1\over\vol_\Gamma}\biggr)^2\muhyp(z)\cr
&={1\over\genus_X^2}\int_{z\in\Gamma\backslash\hyp}F_\Gamma(z)^2\muhyp(z)
-{1\over\vol_\Gamma}\cr
&={1\over\genus_X}\int_{\Gamma\backslash\hyp}F_\Gamma\mucan_X
-{1\over\vol_\Gamma}.}
\eqnumber{eq:Parseval}
$$
Together with Lemma~\ref{lemma:MGamma} and the definition
of~$\zeta_\Gamma$, this finishes the proof.\endproof

We now extend our bounds on~$h_\Gamma$ to the discs around the cusps.

\proclaim Lemma. Let $\Gamma$ be a cofinite Fuchsian group, and let
$X$ be the compactification of\/ $\Gamma\backslash\hyp$.  Let $\fc$ be
a cusp of\/~$\Gamma$, and let $\epsilon$ be a real number
satisfying~\eqref{ineq:eps-sufficiently-small}.  For all $z\in
D_\fc(\epsilon)$, we have
$$
h_{\Gamma,\fc}^-(z)\le h_\Gamma(z)\le h_{\Gamma,\fc}^+(z),
$$
where
$$
h_{\Gamma,\fc}^+(z)=\sup_{\partial\bar D_\fc(\epsilon)} h_\Gamma
+{1\over \vol_\Gamma}\log(\epsilon y_\fc(z))
$$
and
$$
h_{\Gamma,\fc}^-(z)=\inf_{\partial\bar D_\fc(\epsilon)}h_\Gamma
-{\sup_{\partial\bar D_\fc(\epsilon)}F_\Gamma\over\genus_X}
\Bigl({\epsilon\over4\pi}\Bigr)^2
\bigl(1-\exp(4\pi/\epsilon-4\pi y_\fc(z))\bigr)
+{1\over \vol_\Gamma}\log(\epsilon y_\fc(z)).
$$

\label{lemma:extension-cusps-h}

\proof We note that
$$
\Delta h_{\Gamma,\fc}^+(z)=-{1\over\vol_\Gamma}
$$
and
$$
\Delta h_{\Gamma,\fc}^-(z)=
{\sup_{\partial\bar D_\fc(\epsilon)}F_\Gamma\over\genus_X}
(\epsilon y_\fc(z))^2\exp(4\pi/\epsilon-4\pi y_\fc(z))
-{1\over\vol_\Gamma}.
$$
By the non-negativity of~$F_\Gamma$ and
Lemma~\ref{lemma:extension-cusps-F}, this implies
$$
\Delta h_{\Gamma,\fc}^+(z)\le
\Delta h_\Gamma(z)\le
\Delta h_{\Gamma,\fc}^-(z).
$$
Therefore $h_{\Gamma,\fc}^+-h_\Gamma$ and $h_\Gamma-h_{\Gamma,\fc}^-$
are subharmonic functions on~$\bar D_\fc(\epsilon)$.  By the maximum
principle for subharmonic functions, each of these functions assumes
its maximum on the boundary.  The definitions of $h_{\Gamma,\fc}^\pm$
imply that these maxima are non-negative.\endproof

Finally, we prove bounds on the integral
$\int_{\Gamma\backslash\hyp}h_\Gamma\mucan_X$.

\proclaim Lemma. Let $\Gamma$ be a cofinite Fuchsian group, let $X$ be
the compactification of\/ $\Gamma\backslash\hyp$, and let $\eta>0$ be
such that the spectrum of~$-\Delta_\Gamma$ is contained in
$\{0\}\cup[\eta,\infty)$.  Then we have
$$
-\zeta_\Gamma/\eta\le\int_{\Gamma\backslash\hyp}h_\Gamma\mucan_X\le 0.
$$

\label{lemma:bound-int-hGamma}

\proof We use the spectral representations \eqref{eq:FGamma-spectral}
and~\eqref{eq:hGamma-spectral}.  We obtain
$$
\eqalign{
\int_{\Gamma\backslash\hyp}h_\Gamma\mucan_X&=
\int_{z\in\Gamma\backslash\hyp}h_\Gamma(z)\cdot
{1\over\genus_X}F_\Gamma(z)\muhyp(z)\cr
&=-\sum_{j\ge1}{|b_j|^2\over\lambda_j}
-\sum_\fc{1\over 4\pi i}\int_{\Re s=1/2}{|b_\fc(s)|^2\over s(1-s)}.}
$$
We note that the right-hand side is non-positive.  Next, the
assumption that the spectrum of $-\Delta_\Gamma$ is contained in
$\{0\}\cup[\eta,\infty)$ implies
$$
\int_{\Gamma\backslash\hyp}h_\Gamma\mucan_X\ge
-{1\over\eta}\Biggl(\sum_{j\ge1}|b_j|^2+\sum_\fc
{1\over 4\pi i}\int_{\Re s=1/2}|b_\fc(s)|^2ds\Biggr).
$$
Together with~\eqref{eq:Parseval}, this proves the claim.\endproof

\medbreak\noindent
{\it Proof of Theorem~\ref{theorem:bounds-grcan}\/.}\enspace
Part (a) of the theorem follows from the comparison
formula~\eqref{comparison-green}, the bound~\eqref{ineq:grGammaAB}
for~$\gr_\Gamma$, the bound on for~$h_\Gamma$ given by
Lemma~\ref{lemma:bound-hGamma-Y}, and the bound on
$\int_{\Gamma\backslash\hyp}h_\Gamma\mucan_X$ given by
Lemma~\ref{lemma:bound-int-hGamma}.

The proof of parts (b)--(d) is similar.  We first note that, by
Lemmata \ref{lemma:bound-hGamma-Y} and~\ref{lemma:extension-cusps-h},
$$
-S-T(\epsilon_\fc)\le h_\Gamma(z)-{1\over\vol_\Gamma}
\log(\epsilon_\fc y_\fc(z))\le S
\quad\hbox{for all }z\in D_\fc(\epsilon_\fc),
$$
and similarly with $\epsilon_\fc'$ in place of~$\epsilon_\fc$.
Instead of~\eqref{ineq:grGammaAB} we now
invoke \citex{green1}{Proposition~5.5}, which gives bounds for the
function~$\gr_\Gamma$ when one or both variables are near a cusp.  As
in the proof of~(a), it remains to apply the
formula~\eqref{comparison-green} and
Lemma~\ref{lemma:bound-int-hGamma}.\endproof

\section Example: congruence subgroups of~$\SL_2(\Z)$

\label{sec:example}

Let $\Gamma$ be a congruence subgroup of~$\SL_2(\Z)$ such that the
corresponding modular curve~$X$ has positive genus.  Let $n$ be the
level of~$\Gamma$, i.e.\ the minimal positive integer with the
property that $\Gamma$ contains the kernel of the reduction map
$\SL_2(\Z)\to\SL_2(\Z/n\Z)$.

We start by fixing the various parameters.  For the parameter~$\delta$
from Theorem~\ref{theorem:bounds-grcan}, we take
$$
\delta=2.
$$
Selberg's eigenvalue conjecture predicts that all non-zero eigenvalues
of~$-\Delta_\Gamma$ are at least 1/4.  This is at present not known to
be true, but the sharpest known lower bound, due to Kim and
Sarnak~\citex{Kim}{Appendix~2}, allows us to take
$$
\eta=975/4096.
$$

We define
$$
\epsilon=(\delta+\sqrt{\delta^2-1})^{-3/2}\approx 0.139
\quad\hbox{and}\quad
\epsilon'=(\delta+\sqrt{\delta^2-1})^{-1/2}\approx 0.518.
$$
Let $Y_0$ denote the compact subset of $\SL_2(\Z)\backslash\hyp$ which
is the image of the strip
$$
\{x+iy\in\hyp\mid |x|\le 1/2\hbox{ and }
\sqrt{3}/2\le y\le 1/\epsilon\}.
$$
For every cusp $\fc$ of~$\Gamma$, we let $m_\fc$ denote the
ramification index of~$\fc$ over the unique cusp~$\infty$
of~$\SL_2(\Z)$; this equals the index of the corresponding maximal
parabolic subgroups considered modulo $\{\pm1\}$.  
For the parameters $\epsilon_\fc$ and~$\epsilon_\fc'$, we take
$$
\epsilon_\fc=m_\fc\epsilon
\quad\hbox{and}\quad
\epsilon_\fc'=m_\fc\epsilon'.
$$
Using the definition of~$C_\fc(\gamma)$, it is not hard to show that
$$
\min_{\textstyle{\gamma\in\Gamma\atop\gamma\not\in\Gamma_\fc}}
C_\fc(\gamma)\ge m_\fc.
$$
This implies that the parameters $\epsilon_\fc$ and~$\epsilon_\fc'$
satisfy the conditions in Theorem~\ref{theorem:bounds-grcan}.  As in
Theorem~\ref{theorem:bounds-grcan}, let $Y$ be the complement of the
discs $D_\fc(\epsilon_\fc)$.  Then $Y$ is the inverse image of~$Y_0$
in $\Gamma\backslash\hyp$.

We will need an upper bound on the point counting function
$N_\Gamma(z,17)$ for $z\in Y_0$.  It is clear from the definition of
$N_\Gamma(z,U)$ that
$$
\sup_{z\in Y}N_\Gamma(z,U)
\le\sup_{z\in Y_0}N_{\SL_2(\Z)}(z,U).
$$
Using the methods of \citex{green1}{\S4.3}, we have
$$
N_{\SL_2(\Z)}(z,17)\le 226
\quad\hbox{for all }z\in Y_0.
\eqnumber{ineq:NY0}
$$
We next compute suitable $A$ and~$B$
satisfying~\eqref{ineq:grGammaAB}.  For this we use \eqref{ineq:NY0}
and the remaining part of \citex{green1}{\S4.3}, with the same
parameters $\alpha^\pm$, $\beta^\pm$ and~$\gamma^\pm$.  The result is
$$
A=-3.00\cdot 10^4
\quad\hbox{and}\quad
B=1.58\cdot 10^4.
$$
We next find a suitable value of the parameter~$C$.  We use
\citex{green1}{Lemma~2.4}, which says
$$
\Phi_\Gamma(z,\lambda)
\le{\pi\over(2\pi-4)^2}N_{\SL_2(\Z)}(z,17)\lambda
\quad\hbox{for all $z\in\hyp$ and all }\lambda\ge1/4.
$$
The inequality~\eqref{ineq:NY0} implies that we can take
$$
C=137.
$$

We continue with explicit bounds on the canonical $(1,1)$-form.  For
the parameter $a$ from Proposition~\ref{prop:bound-F}, we take
$$
a=1.44.
$$
Again using the method from \citex{green1}{\S4.3}, we compute an upper
bound for $N_{\SL_2(\Z)}(z,2a^2-1)$ for $z\in Y_0$.  The result is
$$
N_{\SL_2(\Z)}(z,2a^2-1)\le 58
\quad\hbox{for all }z\in Y_0.
$$
Substituting this in the bound from Proposition~\ref{prop:bound-F}, we
see that
$$
\sup_Y F_\Gamma\le 25.7.
$$
For every cusp~$\fc$, Lemma~\ref{lemma:extension-cusps-h} implies
$$
\sup_{D_\fc(\epsilon_\fc)}F_\Gamma\le\max\biggl\{1,
\Bigl({\epsilon_\fc\over2\pi}\Bigr)^2\biggl\}\sup_Y F_\Gamma.
$$
From the definition of~$\epsilon_\fc$ and the fact that all
ramification indices~$m_\fc$ are bounded by the level~$n$ of~$\Gamma$,
we conclude
$$
\eqalign{
\sup_X F_\Gamma&\le\max\biggl\{1,
\Bigl({n\epsilon\over2\pi}\Bigr)^2\biggr\}\sup_Y F_\Gamma\cr
&\le\max\{25.7,0.0126n^2\}.}
$$

Finally, we consider the invariant~$\zeta_\Gamma$.  Using
$\int_X\mucan_X=1$ and $\genus_X\ge1$, we make the rather coarse
estimate
$$
\zeta_\Gamma\le\sup_X F_\Gamma\le\max\{25.7,0.0126 n^2\}.
$$

\medbreak\noindent
{\it Proof of Theorem~\ref{theorem:example}\/.}\enspace With the above
estimates, we obtain the following bounds on the various constants in
the theorem:
$$
\eqalign{
S&\le\max\{172,3.79n\},\cr
T(\epsilon_\fc)&\le 0.00313 n^2,\cr
T(\epsilon_\fc')&\le 0.0436 n^2,\cr
\tilde A_\fc&\ge -3.00\cdot 10^4-0.0279n,\cr
\tilde B_\fc&\le 1.58\cdot 10^4+0.0279n.}
$$
The theorem follows from Theorem~\ref{theorem:bounds-grcan} and the
above bounds.\endproof

\unnumberedsection References

\def\PSL{\mathop{\rm PSL}\nolimits}

\normalparindent=25pt
\parindent=\normalparindent
\parskip=1ex plus 0.5ex minus 0.2ex

\reference{Arakelov} {\cyr\S}.~{\cyr\Yu}.~{\cycsc \A\r\a\k\ye\l\o\v},
{\cyr \T\ye\o\r\i\ya\ \p\ye\r\ye\s\ye\ch\ye\n\i\j\ \d\i\v\i\z\o\r\o\v\ \n\a\
\a\r\i\f\m\ye\t\i\ch\ye\s\k\o\j\ \p\o\v\ye\r\kh\n\o\s\t\i}.
{\cyti \I\z\v\ye\s\t\i\ya\ \A\k\a\d\ye\m\i\i\ \N\a\u\k\ \S\S\S\R{\rm
  ,} \s\ye\r\i\ya\ \m\a\t\ye\m\a\t\i\ch\ye\s\k\a\ya\/} {\bf 38}
  (1974), {\cmcyr\char"19}~6, 1179--1192.
\hfil\break
S. Yu.\ {\sc Arakelov}, Intersection theory of divisors on an
arithmetic surface.  {\it Mathematics of the USSR Izvestiya\/} {\bf 8}
(1974), 1167--1180.  (English translation.)

\reference{thesis} P. J. {\sc Bruin}, {\sl Modular curves, Arakelov
theory, algorithmic applications\/}.  Proefschrift (Ph.\thinspace D.\
thesis), Universiteit Leiden, 2010.

\reference{green1} P. J. {\sc Bruin}, Explicit bounds on Green
functions of Fuchsian groups.  Submitted; preprint available at {\tt
http://arxiv.org/abs/1205.6306}.

\reference{book} S. J. {\sc Edixhoven} and J.-M. {\sc Couveignes}
(with R. S. {\sc de Jong}, F. {\sc Merkl} and J. G. {\sc Bosman}),
{\sl Computational Aspects of Modular Forms and Galois
Representations\/}.  Annals of Mathematics Studies {\bf 176}.
Princeton University Press, 2011.



\reference{Faltings: Calculus on arithmetic surfaces} G. {\sc
Faltings}, Calculus on arithmetic surfaces.  {\it Annals of
Mathematics\/}~(2) {\bf 119} (1984), 387--424.

\reference{Iwaniec} H. {\sc Iwaniec}, {\sl Introduction to the
Spectral Theory of Automorphic Forms\/}.  Revista Mate\-m\'a\-tica
Ibero\-americana, Madrid, 1995.

\reference{Fay} J. D. {\sc Fay}, Fourier coefficients of the resolvent
for a Fuchsian group.  {\it Journal f\"ur die reine und angewandte
Mathematik\/} {\bf 293/294} (1977), 143--203.

\reference{Hejhal-I} D. A. {\sc Hejhal}, {\sl The Selberg trace
formula for~$\PSL(2,{\bf R})$\/}, Volume~1.  Lecture Notes in
Mathematics~{\bf 548}.  Springer-Verlag, Berlin/Heidelberg, 1976.

\reference{Hejhal-II} D. A. {\sc Hejhal}, {\sl The Selberg trace
formula for~$\PSL(2,{\bf R})$\/}, Volume~2.  Lecture Notes in
Mathematics~{\bf 1001}.  Springer-Verlag, Berlin/Heidelberg, 1983.

\reference{Jorgenson-Kramer: automorphic forms} J. {\sc Jorgenson} and
J. {\sc Kramer}, Bounding the sup-norm of automorphic forms.  {\it
  Geometric and Functional Analysis\/} {\bf 14} (2005), no.~6,
1267--1277.

\reference{Jorgenson-Kramer: Green's functions} J. {\sc Jorgenson} and
J. {\sc Kramer}, Bounds on canonical Green's functions.  {\it
Compositio Mathematica\/} {\bf 142} (2006), no.~3, 679--700.


\reference{Kim} H. H. {\sc Kim}, Functoriality for the exterior square
of~${\rm GL}_4$ and the symmetric fourth of~${\rm GL}_2$.  With
appendix~1 by D. {\sc Ramakrishnan} and appendix~2 by {\sc Kim} and
P.~{\sc Sarnak}.  {\it Journal of the A.M.S.\/} {\bf 16} (2002),
no.~1, 139--183.





\vskip12mm

\leftline{Peter Bruin}
\leftline{Institut f\"ur Mathematik}
\leftline{Universit\"at Z\"urich}
\leftline{Winterthurerstrasse 190}
\leftline{CH-8057 Z\"urich}
\smallskip
\leftline{\tt peter.bruin@math.uzh.ch}

\bye